\documentclass[11pt, reqno]{amsart}

\usepackage{amscd,amsmath}
\usepackage{mathrsfs}
\usepackage{amsfonts}
\usepackage{amssymb}
\usepackage{enumerate}
\usepackage{setspace}
\usepackage{color}
\usepackage{esint}

\usepackage[margin=1in]{geometry}

\usepackage[english]{babel}

\usepackage{tabu}

\usepackage[colorlinks=true, citecolor=blue]{hyperref}
\newtheorem{theorem}{Theorem}[section]
\newtheorem{lemma}{Lemma}[section]

\newtheorem{definition}{Definition}[section]

\newcommand{\eqn}{\begin{eqnarray}}
\newcommand{\een}{\end{eqnarray}}
\DeclareMathOperator*{\esssup}{ess\,sup}

\numberwithin{equation}{section}

\begin{document}

\title[DLSS-type equations in critical spaces]{Global existence and exponential decay to equilibrium for DLSS-type equations}

\author{Hantaek Bae}
\address{Department of Mathematical Sciences, Ulsan National Institute of Science and Technology (UNIST), Republic of Korea}
\email{hantaek@unist.ac.kr}

\author{Rafael Granero-Belinch\'on}
\address{Departamento de Matem\'aticas, Estad\'istica y Computaci\'on, Universidad de Cantabria, Spain.}
\email{rafael.granero@unican.es}

\date{\today}



\begin{abstract}
In this paper, we deal with two logarithmic fourth order differential equations: the extended one-dimensional DLSS equation and its multi-dimensional analog. We show the global existence of solution in critical spaces, its convergence to equilibrium and the gain of spatial analyticity for these two equations in a unified way.  
\end{abstract}

\maketitle

%{\small \tableofcontents}

\vspace{-2ex}

%%%%%%%%%%%%%%%%%%%%%%%%%
\section{Introduction}
%%%%%%%%%%%%%%%%%%%%%%%%%%
Fourth order differential equations appear in many applications such as thin films (see for instance the works by Constantin, Dupont, Goldstein, Kadanoff, Shelley \& Zhou \cite{constantin1993droplet} and Bertozzi \& Pugh  \cite{Bertozzi1996}), crystal surface models (see the works by Krug, Dobbs \& Majaniemi \cite{krug1995adatom} and Marzuola \& Weare \cite{marzuola2013relaxation}) and quantum semiconductors (see the papers by Ancona \cite{ancona1987diffusion} and Gasser, Markowich, Schmidt \& Unterreiter \cite{gasser1995macroscopic}). In particular, Derrida, Lebowitz, Speer \& Spohn derived the following logarithmic fourth order equation (DLSS in short) as a model of interface fluctuations in a certain spin system \cite{DLSS2}
\eqn \label{DLSS}
w_{t}+\partial_{x}^2(w\partial_{x}^2(\log w))=0.
\een
This equation is a nonlinear parabolic equation. Although the theory of second-order diffusion equations is well known, there are few mathematical results for higher-order equations, and even less results addressing high-order equation in the case of several spatial dimensions. In this paper we first consider the multi-dimensional DLSS equation taking the form \cite{ancona1987diffusion,degond2005quantum,gasser1995macroscopic},
\eqn  \label{hDLSS}
w_{t}+\sum^{d}_{i,j=1}\partial_{i}\partial_{j}(w\partial_{i}\partial_{j}(\log w))=0,
\een 
where the dimension satisfies $d=1,2$ or $d=3$ and the spatial variable $x$ lies in the $d$- dimensional flat torus, or equivalently, on $[-\pi,\pi]^d$ with the periodic boundary conditions. We also study the extended 1D DLSS equation derived by Bordenave, Germain \& Trogdon \cite{Germain}
\eqn \label{eDLSS}
w_{t}-\frac{\mu \Gamma}{3} \partial_{x}^3w-\mu \Gamma \partial_{x}(w \partial_{x}^2\log w))=\epsilon(2\Gamma^{2}-2\Gamma)\partial_{x}^2(w\partial_{x}^2(\log w)),
\een
where $0\leq \Gamma\leq \frac{1}{4}$ and $-1\leq \mu\leq 1$. We note that (\ref{eDLSS}) is reduced to the 1D (\ref{DLSS}) for $\mu=0$, $\epsilon=\frac{8}{3}$, and $\Gamma=\frac{1}{4}$, while (\ref{eDLSS}) becomes \cite[eq.(10)]{DLSS1} for $\epsilon=0$.

Besides the applications to quantum semiconductors and spin systems mentioned above,  \eqref{DLSS} is also intereting due to the fact that, as shown by Gianazza, Savar{\'e} \& Toscani \cite{gianazza2009wasserstein}, it  is a Wasserstein gradient flow for the Fisher information 
$$
\mathcal{F}(u)=-\int \Delta u \log (u) dx
$$
and so \eqref{DLSS} can be written in the form
$$
u_t=\nabla\cdot \left(u \nabla \frac{\delta \mathcal{F}}{\delta u}\right).
$$
This somehow establishes a link between the DLSS equation and the heat equation \cite{jungel2007review}. Indeed, the heat equation is a Wasserstein gradient flow for the entropy
$$
\mathcal{H}(u)=\int u\log(u)dx,
$$
and its entropy production is the widely appearing Fisher information.

The available literature studying \eqref{DLSS} and \eqref{hDLSS} is large. First, Bleher, Lebowitz \& Speer \cite{bleher1994existence} proved that \eqref{DLSS} has local classical solutions starting from $H^{1}(\mathbb{T})$ \emph{positive} initial data. These solutions remain smooth as long as the solution remains positive. We note that, although positivity preservation is a basic property for second order diffusions, it is no longer true for higher order diffusions. In fact, we can see numerical simulations violating the maximum principle in \cite{jungel2007review}. This leaves the door open to possible singularity formation as $u$ approaches zero at some point. There are some papers showing the global existence of \emph{non-negative} weak solutions. These results are mainly based on some appropriate Lyapunov functionals. In particular, both the entropy $\mathcal{H}$ and Fisher information $\mathcal{F}$ decay and can be used to obtain a priori estimates. For \eqref{DLSS}, this approach was exploited by J{\"u}ngel \& Pinnau \cite{jungel2000global} (see also Pia Gualdani, J{\"u}ngel \& Toscani \cite{pia2006nonlinear}) where the authors used the functional
$$
\mathcal{G}(t)=\int u(x,t)-\log(u(x,t))dx.
$$
In the multi-dimensional case \eqref{hDLSS}, the global existence of \emph{non-negative} weak solution was obtained by Gianazza, Savar{\'e} \& Toscani \cite{gianazza2009wasserstein} and J{\"u}ngel \& Matthes \cite{jungel2008derrida}.

In addition to the existence results, there are also some works studying the decay to the equilibrium. In the one-dimensional case, C{\'a}ceres, Carrillo \& Toscani \cite{caceres2005long}, Dolbeault, Gentil \& J{\"u}ngel \cite{dolbeault2006logarithmic}, J{\"u}ngel \& Violet \cite{jungel2007first} and J{\"u}ngel \& Matthes \cite{jungel2008derrida} proved the decay of appropriate functionals in the case of periodic conditions. The case of Dirichlet and Neumann boundary conditions in one-dimension has been studied  by J{\"u}ngel \& Pinnau \cite{jungel2000global} and J{\"u}ngel \& Toscani \cite{jungel2003exponential}. In the multi-dimensional setting, the reader is referred to \cite{gianazza2009wasserstein,jungel2008derrida}. To the best of our knowledge, the only decay result in the multi-dimensional setting with periodic boundary conditions is \cite{jungel2008derrida}, where the authors proved decay of $L^p$-norm-like functionals for 
\begin{equation}\label{decay}
1\leq p<\frac{(\sqrt{d}+1)^2}{d+2}<\infty.
\end{equation}

The uniqueness question has also been studied. Fischer \cite{fischer2013uniqueness} proved that weak solutions in a certain class are unique. Remarkably, this class contains the weak solutions constructed in \cite{jungel2008derrida}. We would like to remark that the question of uniqueness is rather subtle. As it was pointed out in \cite{jungel2007review}, there are non-negative explicit functions that are steady solutions of \eqref{hDLSS}. These non-negative functions lead, after invoking the previously mentioned existence results, to a time-dependent weak solution converging to the homogeneous steady state. In other words, starting from smooth initial data we can have two different weak solutions. 

The results for equation \eqref{eDLSS} are more scarce. This equation was derived and studied by Bordenave, Germain \& Trogdon \cite{Germain}. In particular, the authors found a number of Lyapunov functionals and used them to study the asymptotic behavior. 

Some other related works are \cite{bukal1,bukal2,bukal3} (and the references therein) where the DLSS equation and some generalization are studied.

In this paper, we deal with (\ref{hDLSS}) and (\ref{eDLSS}) in a unified way. Assuming that the initial data satisfy certain \emph{explicit} size restrictions in the Wiener algebra $\mathcal{A}_{0}$ (defined in \eqref{Wiener space}), we prove the global existence of solutions, its convergence to the steady state in $L^{\infty}$ and instantaneous gain of analyticity. We will present our results in Section \ref{sec 2}. Then, we prove Theorem \ref{DLSS GW 1} and \ref{DLSS GW 2} in Sections \ref{sec 4} and \ref{sec 5}, respectively. 

%Then, in Appendix \ref{A: alternative analyticity} we proved another proof of the gain of analyticity.

We also want to mention that our approach has the following advantages:

\begin{itemize}
\item We would like to emphasize that the smallness conditions of initial data are given \emph{explicitly} in terms of the parameters in the equation.
\item If $w$ is a solution to \eqref{hDLSS}, then the rescaled function $w_\lambda(t,x)=w(\lambda^4 t, \lambda x)$ is another solution for every $\lambda>0$.  Some spaces with scaling-invariant norm are 
$$
L^\infty, \;H^{d/2},\;\mathcal{A}_{0},\;\text{etc}.
$$
Thus, our results involve a critical space for \eqref{hDLSS}.  We should compare our global existence results with the one contained in \cite{bleher1994existence} where the size constrain is at the level of $\dot{H}^1(\mathbb{T}^1)$. In particular, the initial data that we consider can be arbitrarily large in the space $\dot{H}^{1}$ while still leading to a global solution. Moreover, to the best of our knowledge, our results are the first results showing analyticity (not merely $C^\infty$ as in \cite{bleher1994existence}) of solutions of (\ref{hDLSS}) and (\ref{eDLSS}). 
\item We prove the exponential decay to equilibrium in the $L^\infty$ norm which generalizes the result in \cite{jungel2008derrida}. In particular, the decay \eqref{decay} is extended to $1\leq p\leq\infty.$
\item Our approach is very flexible and can be implemented to other high order semi-linear or quasi-linear equations  in an arbitrary dimension $d$. We refer the  reader to \cite{gan2,gan3,gan1} for a free boundary problem arising in the dynamics of a fluid in a porous medium, to \cite{granero2018asymptotic} for  a nonlocal quasilinear diffusion, to \cite{burczak2016generalized} for  the doubly parabolic Keller-Segel system, to \cite{bruell2019thin} for thin film equations, and to \cite{ambrose2018radius,granero2019global,liu2019global} for  the evolution of crystal surfaces.  These results are obtained by mainly observing that the Wiener algebra norm is a Lyapunov functional regardless of the order of the diffusion, the local/nonlocal character of nonlinear terms and the dimension $d$. 
\end{itemize}

%%%%%%%%%%%%%%%%%
\section{Main results}\label{sec 2}
%%%%%%%%%%%%%%%%%

%%%%%%%%%%%%%%%%%%%%%%%%%%%%
\subsection{Notation and definitions}
%%%%%%%%%%%%%%%%%%%%%%%%%%%%
Before stating our results, we fix some notation and introduce the functional spaces that we will use.

The spatial derivatives are denoted by 
$$
\partial_{i}f=f,_{i},
$$ 
$$
\partial_{i}\partial_{j}f=f,_{ij},
$$
and so on, where we also used Einstein convention for the repeated indices. In particular, $$
\Delta f=f,_{ii}
$$ and 
$$
\Delta^{2}f=f,_{ijij}.
$$ 
In the one dimensional case, we also write 
$$
f_x=f,_1.
$$ 
Similarly, the time derivative is denoted as 
$$
f_t=\partial_{t}f.
$$ 

For  $n\in \mathbb{N}$ we denote by
\begin{equation*}
W^{n,p}=\left\{u\in L^p\text{ such that } \|u\|_{L^p}^p+\|\partial_x^n u\|_{L^p}^p<\infty
\right\}
\end{equation*}
the standard $L^p$-based Sobolev spaces. Then we define the norm as
$$
\|u\|_{W^{n,p}}^p:= \|u\|_{L^p}^p+\|\partial_x^n u\|_{L^p}^p. 
$$
When $p=2$, we use $W^{n,2}=H^{n}$. 

The $k$-th Fourier coefficients of a $2\pi$-periodic function on $\mathbb{T}^{d}$ are 
\[
\widehat{u}(k)=\frac{1}{(2\pi)^{d}}\int_{\mathbb{T}^{d}}u(x)e^{-ik\cdot x}dx,
\]
and the Fourier series expansion of $u$ is given by
\[
u(x)=\sum_{k\in \mathbb{Z}^{d}}\widehat{u}(k) e^{ik\cdot x}.
\]
Using this, we define the Wiener spaces: for $s\ge 0$
\eqn \label{Wiener space}
\mathcal{A}_{s}=\left\{u\in L^{1}(\mathbb{T}^{d})\text{ such that } \sum_{k\in \mathbb{Z}^{d}}|k|^{s}\left|\widehat{u}(k)\right|<\infty\right\}. 
\een
Then we define the following norm in this space
$$
|u|_{s}:=\sum_{k\in \mathbb{Z}^{d}}|k|^{s}\left|\widehat{u}(k)\right|
$$
We note that  $\mathcal{A}_{0}$ is a Banach algebra. Moreover, 
$$
\mathcal{A}_{s}\subset C^s\subset H^{s}.
$$ 
Being $\partial$ a differential operator of order 1 we also have the following properties  
\begin{equation}
\begin{split}
& \left|\partial^{l}u\right|_{0}\leq |u|_{l},\\
& \left|\partial^{l}u \partial^{l'}u\right|_{0}\leq |u|_{l} |u|_{l'},\\
& |u|_{l} \leq |u|_{l'} \quad \text{if $\ l'\ge l>0$}.
\end{split}
\end{equation}
We finally have the following interpolation relationships:
\eqn \label{interpolation}
|u|_{s}\leq |u|^{1-\theta}_{0} |u|^{\theta}_{r} \quad \text{for all $0\leq s\leq r$ with} \  \theta=\frac{s}{r}.
\een

Let $X$ be a Banach space. Then, $L^{p}_{T}X$ denotes the Banach set of Bochner measurable functions $f$ from $(0,T)$ to $X$ such that 
$$
\|f(t)\|_{X}\in L^{p}(0,T)
$$ 
endowed with the norm 
\begin{equation*}
\begin{split}
\left(\int ^{T}_{0} \|f(\cdot, t)\|^{p}_{X}dt\right)^{\frac{1}{p}} \ \text{for $1\leq p<\infty$ \ \ },
\end{split}
\end{equation*}
or
\begin{equation*}
\begin{split}
\esssup_{0\leq t\leq T} \|f(\cdot,t)\|_{X} \quad \text{for $p=\infty$}.
\end{split}
\end{equation*}

%%%%%%%%%%%%
\subsection{Setting}
%%%%%%%%%%%
We assume that the spatial variable $x$ lies in the $d$- dimensional flat torus $[-\pi,\pi]^d$ with the periodic boundary conditions.  Together with \eqref{hDLSS} and \eqref{eDLSS}, we have to consider \emph{non-negative} initial data $w(x,0)=w_{0}(x)$. We note that both (\ref{hDLSS}) and (\ref{eDLSS}) preserve the mean $\langle w \rangle$:
\[
\langle w(t) \rangle=\int_{\mathbb{T}^{d}}w(t,x)dx=\int_{\mathbb{T}^{d}}w_{0}(x)dx=\langle w_{0} \rangle>0
\]
In this setting, this mean will then correspond to the steady state. Without loss of generality, we assume $\langle w_{0} \rangle=1$.

%%%%%%%%%%%%%%%%%%%%%%%%%%%%
\subsection{The multi-dimensional DLSS equation}
%%%%%%%%%%%%%%%%%%%%%%%%%%%

We now reformulate (\ref{hDLSS}) in a new variable which is more convenient to show our results. Let $u=w-1$. Then, $u$ satisfies 
\eqn \label{eq of u}
u_{t} +((u+1)(\log (u+1)),_{ij}),_{ij}=0,\quad 1\leq i,j\leq d,
\een
with initial data $u(x,0)=u_0(x)=w_0(x)-1$. We note that 
\[
\langle u(t) \rangle=\int_{\mathbb{T}^{d}}u(t,x)dx=\int_{\mathbb{T}^{d}}u_{0}(x)dx=\langle u_{0} \rangle=0.
\]
This property will be used several times in proving our results.

\begin{definition}
We say that $u \in L^{\infty}_{T}L^{\infty}\cap L^{2}_T H^{1}\cap L^1_{T}W^{2,1}$  is a strong solution of (\ref{eq of u}) if 
\begin{equation}\label{weak form}
\begin{split}
&-\int^{T}_{0}\int_{\mathbb{T}^{d}} \partial_{t}\phi(t,x)u(t,x)dxdt-\int_{\mathbb{T}^{d}}\phi(0,x)u_{0}(x)dx\\
&+\int^{T}_{0}\int_{\mathbb{T}^{d}}\log(u(t,x)+1),_{ij}(u(t,x)+1)\phi,_{ij}(t,x)dxdt =0
\end{split}
\end{equation}
for all $\phi \in C^{\infty}([0,\infty)\times \mathbb{T}^{d})$.
\end{definition}

We note 
\[
\log(u+1),_{ij}(u+1)=u_{,ij}-\frac{u_{,i} u_{,j}}{1+u}.
\]
From this, the regularity condition $u \in L^{\infty}_{T}L^{\infty} \cap L^{2}_T H^{1} \cap L^1_{T}W^{2,1}$ seems  \emph{minimal} to define (\ref{weak form})  by the following reasons: 
\begin{itemize}
\item $u \in L^1_{T}W^{2,1}$ is required to define $u_{,ij}$ weakly, 
\item $u \in L^{\infty}_{T}L^{\infty}$, which will be sufficiently small, to make $1+u>0$; \item $u \in L^{2}_{T}H^{1}$ to define $u_{,i} u_{,j}$ weakly. 
\end{itemize}
This required regularity of $u$ will be achieved by taking a sufficiently small $u_{0}\in \mathcal{A}_{0}$. The smallness condition of $u$ is defined explicitly by  the following the rational function:
\eqn \label{P1}
P_{1}(z)=\frac{4z}{1- z}+\frac{5z^2}{(1- z)^{2}}+\frac{4z^{3}}{(1- z)^{3}}.
\een

\begin{theorem} \label{DLSS GW 1}
Let $u_{0} \in \mathcal{A}_{0}$ be a zero mean function that satisfies $|u_{0}|_{0}<1$ and consider $P_1(z)$ as in \eqref{P1}. Assume that the initial data satisfies 
$$
d^{2}P_{1}(|u_{0}|_{0})<1.
$$ 
Then there is a solution 
$$
u\in L^\infty([0,T];L^\infty)\cap L^{p}(0,T;\mathcal{A}_{\frac{4}{p}}),\;\;1<p<\infty,
$$
of (\ref{eq of u})  for all $T>0$. Moreover, $u$ converges uniformly to $0$ exponentially in time:
\[
\|u(t)\|_{L^\infty}\leq |u_{0}|_{0}e^{-(1-d^{2}P_{1}(|u_{0}|_{0}))t}.
\]
Finally, $u$ becomes instantaneously spatial analytic with increasing radius of analyticity, namely,
\[
e^{\sigma t|k|}\widehat{u}(t,k)\in L^\infty(0,T;\ell^1)
\]
for a sufficiently small $\sigma>0$ such that $1-\sigma-d^{2}P_{1}(|u_{0}|_{0})>0$. 
\end{theorem}

%%%%%%%%%%%%%%%%%%%%%%%%
\subsection{The extended 1D DLSS equation}
%%%%%%%%%%%%%%%%%%%%%%%
As before, we define $u=w-1$. Then, $u$ satisfies the following equation:
\begin{align} \label{u ext}
u_{t}-\frac{\mu\Gamma}{3} u_{xxx}-\mu\Gamma ((u+1)(\log (u+1))_{xx})_{x}=-\epsilon(2\Gamma-2\Gamma^{2})((u+1)(\log (u+1))_{xx})_{xx}.
\end{align}

\begin{definition}
We say that $u \in L^{\infty}_{T}L^{\infty}\cap L^{2}_T H^{1} \cap L^1_{T}W^{2,1}$  is a strong solution of  (\ref{u ext})  if 
\begin{equation}
\begin{split}
&-\int^{T}_{0}\int_{\mathbb{T}}\partial_{t}\phi(t,x)u(t,x)dxdt+\frac{\mu \Gamma}{3}\int^{T}_{0}\int_{\mathbb{T}}u(t,x)\phi_{xxx}(t,x)dxdt -\int_{\mathbb{T}}\phi(0,x)u_{0}(x)dx\\
&+\mu \Gamma \int^{T}_{0}\int_{\mathbb{T}}(\log (u(t,x)+1))_{xx}(1+u(t,x))\phi_x(t,x)dxdt\\
&+\epsilon(2\Gamma^{2}-2\Gamma)\int^{T}_{0}\int_{\mathbb{T}}(\log (u(t,x)+1))_{xx}(1+u(t,x))\phi_{xx}(t,x)dxdt=0
\end{split}
\end{equation}
for all $\phi \in C^{\infty}([0,\infty)\times \mathbb{T}^{d})$.
\end{definition}

We then follow the arguments used for (\ref{eq of u}) to state the following result. Let
\eqn \label{P2}
P_2(z)=\frac{z}{1- z}+\frac{z^2}{(1- z)^2}.
\een

\begin{theorem} \label{DLSS GW 2}
Let $u_{0} \in \mathcal{A}_{0}$ be a zero mean function that satisfies $|u_{0}|_{0}<1$ and consider $P_j(z)$, $j=1,2$, as in \eqref{P1} and \eqref{P2}. Assume that the initial data satisfies 
$$
\mu \Gamma P_{2}(|u_{0}|_{0})+\epsilon(2\Gamma^{2}-2\Gamma) P_{1}(|u_{0}|_{0})<\epsilon(2\Gamma^{2}-2\Gamma).
$$
Then there is a solution  
$$
u\in L^\infty([0,T];L^\infty)\cap L^{p}(0,T;\mathcal{A}_{\frac{4}{p}}),\;1<p<\infty,
$$
of (\ref{u ext}) for all $T>0$. Moreover, $u$ converges uniformly to $0$ exponentially in time
\[
\|u(t)\|_{L^\infty}\leq |u_{0}|_{0}e^{-(\epsilon(2\Gamma^{2}-2\Gamma)-\mu \Gamma P_{2}(|u_{0}|_{0})-\epsilon(2\Gamma^{2}-2\Gamma) P_{1}(|u_{0}|_{0}))t}.
\]
Finally, $u$ becomes instantaneously spatial analytic with increasing radius of analyticity, namely
\[
e^{\sigma t|k|}\widehat{u}(t,k)\in L^\infty(0,T;\ell^1)
\]
for a  sufficiently small $\sigma>0$ such that $\epsilon(2\Gamma^{2}-2\Gamma)-\sigma-\mu \Gamma P_{2}(|u_{0}|_{0})-\epsilon(2\Gamma^{2}-2\Gamma)P_{1}(|u_{0}|_{0})>0$. 
\end{theorem}

%%%%%%%%%%%%%%%%%%%%%%%%%%%%%%%%%%%%%%%%
\section{Proof of Theorem \ref{DLSS GW 1}}\label{sec 4}
%%%%%%%%%%%%%%%%%%%%%%%%%%%%%%%%%%%%%%%%

%%%%%%%%%%%%%%%%%%%%%%%%
\subsection{A priori estimates}\label{apriori}
%%%%%%%%%%%%%%%%%%% %%%%
We compute
\begin{align*}
u_{t}&=-((u+1)(\log (u+1)),_{ij}),_{ij}\\
&=-\left((u+1)\left(\frac{u,_{ij}}{1+u}-\frac{u,_ju,_{i}}{(1+u)^2}\right)\right),_{ij}\\
&=-u,_{ijij}+\left(\frac{u,_j u,_{i}}{1+u}\right),_{ij}\\
&=-u,_{ijij}+\frac{u,_{jj} u,_{ii}+u,_j u,_{iij}+u,_{ijj} u,_{i}+u,_{ij} u,_{ij}}{1+u}\\
&\quad
-\frac{u,_{ii} u,_{j}u,_{j}+3u,_{ij} u,_{i}u,_{j}+u,_{jj} u,_{i}u,_{i}}{(1+u)^2}+2\frac{u,_{j}u,_{j} u,_{i}u,_{i}}{(1+u)^3}.
\end{align*} 
Thus, we can rewrite the equation as 
\begin{equation}\label{eq:urevised}
\begin{split}
u_{t}+\Delta^{2}u&=\text{I}_1+\text{I}_2+\text{I}_3,
\end{split}
\end{equation}
where
$$
\text{I}_1=\frac{u,_{jj} u,_{ii}+u,_j u,_{iij}+u,_{ijj} u,_{i}+u,_{ij} u,_{ij}}{1+u},
$$
$$
\text{I}_2=-\frac{u,_{ii} u,_{j}u,_{j}+3u,_{ij} u,_{i}u,_{j}+u,_{jj} u,_{i}u,_{i}}{(1+u)^2},
$$
$$
\text{I}_3=\frac{2u,_{j}u,_{j} u,_{i}u,_{i}}{(1+u)^3}.
$$
By the hypothesis $|u_0|_0<1$, we take the Taylor expansion of the rational functions in $\text{I}_1$, $\text{I}_2$, and $\text{I}_3$:
\begin{equation*}\label{NLs}
\begin{split}
\text{I}_1&=\left(u,_{jj} u,_{ii}+u,_j u,_{iij}+u,_{ijj} u,_{i}+u,_{ij} u,_{ij}\right)\sum^{\infty}_{n=0}(-1)^{n}u^{n}\\
\text{I}_2&=\left(u,_{ii} u,_{j}u,_{j}+3u,_{ij} u,_{i}u,_{j}+u,_{jj} u,_{i}u,_{i}\right)\sum^{\infty}_{n=1}(-1)^{n}nu^{n-1}\\
\text{I}_3&=2(u,_{j}u,_{j} u,_{i}u,_{i})\sum^{\infty}_{n=2}(-1)^{n}n(n-1)u^{n-2}.
\end{split}
\end{equation*}
Writing
$$
\mathcal{S}_{1}=\sum^{\infty}_{n=0}(-1)^{n}u^{n},
$$
$$
\mathcal{S}_{2}=\sum^{\infty}_{n=1}(-1)^{n}nu^{n-1},
$$
$$
\mathcal{S}_{3}=\sum^{\infty}_{n=2}(-1)^{n}n(n-1)u^{n-2},
$$
and using the convolution theorem for the Fourier series together with Tonelli's theorem for exchanging the order of summation, we can obtain the following estimates:
\begin{equation}\label{geometric series}
\begin{split}
|\mathcal{S}_{1}|_{0}=&\sum _{k\in \mathbb{Z}^{d}}\left|\sum^{\infty}_{n=0}(-1)^{n}\widehat{u^{n}}(k)\right|\leq \sum^{\infty}_{n=0} |u|^{n}_{0} =\frac{1}{1- |u|_{0}}\\
|\mathcal{S}_{2}|_{0}=&\sum _{k\in \mathbb{Z}^{d}}\left|\sum^{\infty}_{n=1}n(-1)^{n}\widehat{u^{n-1}}(k)\right|\leq \sum^{\infty}_{n=1} n|u|^{n-1}_{0}=\frac{1}{(1- |u|_{0})^{2}}\\
|\mathcal{S}_{3}|_{0}=&\sum _{k\in \mathbb{Z}^{d}}\left|\sum^{\infty}_{n=2}n(n-1)(-1)^{n}\widehat{u^{n-2}}(k)\right|\leq \sum^{\infty}_{n=2} n(n-1)|u|^{n-2}_{0} =\frac{2}{(1- |u|_{0})^{3}}.
\end{split}
\end{equation}
Moreover, for fixed $i,j$, we apply (\ref{interpolation}) to derive 
\begin{equation}\label{DLSS interpolation}
\begin{split}
\left|u,_{i}u,_{ijj}\right|_{0}\leq \left|u\right|_1\left|u\right|_{3}&\leq |u|_{0}|u|_{4},\\
\left|u,_{ij}u,_{ij}\right|_{0}\leq \left|u\right|_2^2& \leq |u|_{0} |u|_{4},\\
\left|u,_{j}u,_{j}u,_{ii}\right|_{0}\leq \left|u\right|_1^2\left|u\right|_{2}& \leq |u|^{2}_{0} |u|_{4},\\
\left|u,_{i}u,_{j}u,_{ij}\right|_{0}\leq \left|u\right|_1^2\left|u\right|_2&\leq |u|^{2}_{0}|u|_{4},\\
\left|u,_{i}u,_{j}u,_{i}u,_{j}\right|_{0}\leq |u|_1^4&\leq |u|^{3}_{0}|u|_{4}.
\end{split}
\end{equation}
By \eqref{eq:urevised}, (\ref{geometric series}) and (\ref{DLSS interpolation}) and summing up in $i,j$, we have 
\eqn \label{DLSS final}
\frac{d}{dt}|u(t)|_{0}+|u(t)|_{4} \leq d^{2}P_{1}(|u(t)|_{0})|u(t)|_{4},
\een
where $P_{1}$ is defined in (\ref{P1}). Suppose $u_{0}$ satisfies 
\eqn \label{DLSS size of u0}
|u_{0}|_{0}<1 \quad \text{and} \quad d^{2}P_{1}(|u_{0}|_{0})<1.
\een
Then, using the monotonicity of $P_1$, we obtain the following a priori bound
\eqn \label{DLSS uniform bound}
|u(t)|_{0}+ (1-d^{2}P_{1}(|u_{0}|_{0}))\int^{t}_{0}|u(s)|_{4}ds \leq |u_{0}|_{0} \quad \text{for all $t\ge 0$.}
\een

%%%%%%%%%%%%%%%%%%%%%
\subsection{Approximate sequence of solutions}
%%%%%%%%%%%%%%%%%%%%%
We  consider the following approximate equation:
\begin{subequations} \label{DLSS N}
\begin{align}
& \partial_{t}u+\Delta^{2}u=-\mathcal{P}_N\left(\left(u_{,i}u_{,j}\sum^{N}_{k=0}(-1)^{k}u^{k}\right)_{,ij}\right), \label{DLSS N a}\\
& u(0,x)=\mathcal{P}_Nu_{0},
\end{align}
\end{subequations}
where $\mathcal{P}_Ng$ is defined as 
\[
\mathcal{P}_Ng=\sum_{k=-N}^N\widehat{g}(k)e^{ikx}.
\]
Since the right-hand side of (\ref{DLSS N a}) is a polynomial of $u$, we can solve (\ref{DLSS N}) using Picard's iteration for the ODE's to obtain a solution $u^{N}$ for each $N$. Following the computations used to derive (\ref{DLSS uniform bound}), we can show that $u^{N}$ satisfies the same estimates as above:
\eqn \label{DLSS Ns}
\frac{d}{dt}\left|u^{N}(t)\right|_{0}+\left|u^{N}(t)\right|_{4} \leq d^{2}P_{1}\left(|u^{N}(t)|_{0}\right)\left|u^{N}(t)\right|_{4}, 
\een
$$
\left|\mathcal{P}_Nu_{0}\right|_{0}\leq |u_{0}|_{0}.
$$
We note that $P_{1}(z)$ is increasing for $0<z<1$. Using this fact with $u_{0}$ satisfying (\ref{DLSS size of u0}),  we obtain 
\eqn \label{DLSS uniform bound N}
\left|u^{N}(t)\right|_{0}+ (1-d^{2}P_{1}(|u_{0}|_{0}))\int^{t}_{0}\left|u^{N}(s)\right|_{4}ds \leq |u_{0}|_{0}
\een
for all $t\ge0$ uniformly in $N\in \mathbb{N}$. In particular we find that
$$
\sup_{N}\left|u^{N}(t)\right|_{0}<1.
$$

%%%%%%%%%%%%%%%%%%%%%%%%%%%%%%%%%%%%%%%%%%%%%%%%%
\subsection{Passing to the limit in the weak formulation}
%%%%%%%%%%%%%%%%%%%%%%%%%%%%%%%%%%%%%%%%%%%%%%%%%

The inequality (\ref{DLSS uniform bound N}) implies that $\{u^{N}\}$ is uniformly bounded in $L^{\infty}_{T}\mathcal{A}_{0} \cap L^{1}_{T}\mathcal{A}_{4}$. By interpolating $L^{\infty}_{T}\mathcal{A}_{0}$ and $L^{1}_{T}\mathcal{A}_{4}$, $\{u^{N}\}$ is also  uniformly bounded in 
\eqn \label{uniform bound of uN 2}
\{u^{N}\}\in L^{p}_{T}\mathcal{A}_{\frac{4}{p}}, \quad 1<p<\infty.
\een 
We observe that (see \cite{cengiz1992duals}), for $1<p\leq \infty$
\begin{equation}\label{dualitywiener}
L^{p}_{T}\ell^1=\left(L^{q}_{T}c_0\right)^*,	\;\;q^{-1}+p^{-1}=1,
\end{equation}
where $c_0$ is the space formed by the sequences whose limit is zero and $\ell^1$ is the space formed by the summable sequences. Due to this, the space
$$
L^{p}_{T}\mathcal{A}_0
$$
is a dual space.

Furthermore, since
\[
\left\|u^{N}\right\|^{2}_{H^{2}}=\sum_{k\in \mathbb{Z}^{d}}|\widehat{u^{N}}(k)|^{2}+ \sum_{k\in \mathbb{Z}^{d}}|k|^{4}|\widehat{u^{N}}(k)|^{2}\leq \left|u^{N}\right|^{2}_{0}+ \left|u^{N}\right|_{0}\left|u^{N}\right|_{4},
\]
$\{u^{N}\}$ is uniformly bounded in 
$$
\{u^{N}\}\in L^{2}_{T}H^{2}.
$$
We also note that  $\left\{\partial_{t}u^{N}\right\}$ is uniformly bounded in $L^{2}_{T}H^{-2}$. Indeed, since $u^{N},_{i}\in L^{4}\mathcal{A}_{0}$ by (\ref{uniform bound of uN 2}), we have
\[
u_{,i}u_{,j}\sum^{N}_{k=0}(-1)^{k}u^{k}\in L^{2}_{T}\mathcal{A}_{0}\subset L^{2}_{T}L^{2}
\]
and thus
\[
\partial_{t}u^{N}=-\Delta^{2}u-\mathcal{P}_N\left(\left(u_{,i}u_{,j}\sum^{N}_{k=0}(-1)^{k}u^{k}\right)_{,ij}\right) \in L^{2}_{T}H^{-2}. 
\]
We now recall Simon's compactness lemma.

\begin{lemma} \cite{simon1986compact} \label{compact lemma}
Let $X_{0}$, $X_{1}$, and $X_{2}$ be Banach spaces such that $X_{0}$ is compactly embedded in $X_{1}$ and $X_{1}$ is a subset of $X_{2}$. Then, for $1\leq p<\infty$, any bounded subset of $\left\{v\in L^{p}_{T}X_{0}: \ \frac{\partial v}{\partial t}\in L^{1}_{T}X_{2}\right\}$ is precompact in $L^{p}_{T}X_{1}$.
\end{lemma}

Then, by Lemma \ref{compact lemma} and Banach-Alaoglu theorem, we have that up to subsequences, still denoted by $\{u^{N}\}$, 
\begin{subequations} \label{DLSS conv 1}
\begin{align}
&u^{N}  \overset{\ast}{\rightharpoonup} u \in L^{\infty}([0,T]\times \mathbb{T}^d), \label{DLSS conv 11}\\
&u^{N}  \overset{\ast}{\rightharpoonup} u \in L^{p}_{T}\mathcal{A}_{\frac{4}{p}}, \quad 1<p<\infty,  \label{DLSS conv 1 b}\\
& u^{N} \rightarrow u \in L^{2}_{T}H^{r}, \quad 0\leq r<2. \label{DLSS conv 1 c},\\
&u^{N}  \overset{\ast}{\rightharpoonup} u \in L^{\infty}([0,T],\mathcal{A}_{0}), \label{DLSS conv 11B}
\end{align}
\end{subequations}

We now rewrite  (\ref{DLSS N}) as 
\[
\partial_{t}u^{N}=\mathcal{P}_N\left(u^{N}_{,ij}-\frac{u^{N}_{,i} u^{N}_{,j}}{1+u^{N}}+\frac{u^{N}_{,i} u^{N}_{,j}}{1+u^{N}}-u_{,i}u_{,j}\sum^{N}_{k=0}(-1)^{k}u^{k}\right)_{,ij}.
\]
We multiply this equation with $\phi$. After integrating by parts and using the Taylor series for the nonlinear term, we obtain 
\begin{equation} \label{weak form N}
\begin{split}
&-\int^{T}_{0}\int_{\mathbb{T}^{d}} \partial_{t}\phi(t,x)u^N(t,x)dxdt-\int_{\mathbb{T}^{d}}\phi(0,x)\mathcal{P}_Nu_{0}(x)dx\\
&+\int^{T}_{0}\int_{\mathbb{T}^{d}}\mathcal{P}_N\left(u^{N}_{,ij}-\frac{u^{N}_{,i} u^{N}_{,j}}{1+u^{N}}\right)\phi,_{ij}(t,x)dxdt\\
&-\int^{T}_{0}\int_{\mathbb{T}^{d}}\mathcal{P}_N\left(u^{N}_{,i} u^{N}_{,j}\sum^{\infty}_{k=N+1}(-u^{N})^{k}\right)\phi,_{ij}(t,x)dxdt =0.
\end{split}
\end{equation}

Due to the previous convergences, we have that
\begin{enumerate}[]
\item \textbullet \ \ By (\ref{DLSS conv 1 b}) with $\widetilde{p}=2$, $u^{N}_{,ij} \rightharpoonup u_{,ij}$ in $L^{2}_{T}L^{2}$. 
\item \textbullet \ \ By (\ref{DLSS conv 1 c}) with $r=\frac{3}{2}$, we have $u^{N}_{,i} \rightarrow u_{,i} \in L^{2}_{T}H^{r-1}\subset L^{2}_{T}L^{q_{1}}$ with $q_{1}=\frac{2d}{d-\frac{1}{2}}$. 
\item \textbullet \ \ By (\ref{DLSS conv 1 b}) with $p=4$, we also have $u^{N}_{,j} \rightharpoonup u_{,j} \in L^{4}_{T}\mathcal{A}_{0}\subset L^{4}_{T}L^{q_{2}}$ with $1\leq q_{2}\leq \infty$. 
\item \textbullet \ \ We finally have $u \in L^{4}_{T}H^{\frac{1}{4}}\subset L^{4}_{T}L^{q_{3}}$ with $q_{3}=\frac{4d}{2d-1}$ which can be obtained by (\ref{DLSS conv 1 c}) with $r=\frac{1}{2}$ and  then interpolated with $L^{\infty}_{T}\mathcal{A}_{0}\subset L^{\infty}_{T}L^{2}$. This implies 
\[
\frac{1}{1+u^{N}}-\frac{1}{1+u}=\frac{u-u^{N}}{(1+u^{N})(1+u)}\rightarrow 0 \quad \text{in} \ \   L^{4}_{T}L^{q_{3}}
\]
because the denominator $(1+u^{N})(1+u)$ does not vanish. 
\end{enumerate}
We now choose $q_{2}$ such that 
\[
\frac{1}{q_{1}}+\frac{1}{q_{2}}+\frac{1}{q_{3}}=1.
\]
Then, we have 
\begin{equation*}
\begin{split}
\int^{T}_{0}\int_{\mathbb{T}^{d}}\mathcal{P}_N\left(u^{N}_{,ij}-\frac{u^{N}_{,i} u^{N}_{,j}}{1+u^{N}}\right)\phi,_{ij}(t,x)dxdt &\rightarrow \int^{T}_{0}\int_{\mathbb{T}^{d}}\left(u_{,ij}-\frac{u_{,i} u_{,j}}{1+u}\right)\phi,_{ij}(t,x)dxdt\\
&=\int^{T}_{0}\int_{\mathbb{T}^{d}}\log(u(t,x)+1),_{ij}(u(t,x)+1)\phi,_{ij}(t,x)dxdt.
\end{split}
\end{equation*}
We finally note that 
\[
\left\|u^{N}_{,i} u^{N}_{,j}\sum^{\infty}_{k=N+1}(u^{N})^{k}\right\|_{L^{2}_{T}L^{2}}\leq  \left\|u^{N}\right\|_{L^{4}_{T}\mathcal{A}_{1}}^2\sum^{\infty}_{k=N+1}\left\|u^{N}\right\|^{k}_{L^{\infty}_{T}\mathcal{A}_{0}} \rightarrow 0 \quad \text{as} \quad N\rightarrow \infty.
\]
We now can pass to the limit to (\ref{weak form N}) to derive (\ref{weak form}) for any $\phi \in C^{\infty}([0,\infty)\times \mathbb{T}^{d})$.

%%%%%%%%%%%%%%%%%%%%%%%%%
\subsection{Decay to equilibrium}
%%%%%%%%%%%%%%%%%%%%%%%%%
Since $\langle u^{N}(t) \rangle=0$, we have $|u^N(t)|_{0}\leq |u^N(t)|_{4}$. Using  (\ref{DLSS Ns}), we obtain
\[
|u^N(t)|_{0}\leq |u_{0}|_{0}e^{-(1-d^{2}P_{1}(|u^N_{0}|_{0}))t}\leq |u_{0}|_{0}e^{-(1-d^{2}P_{1}(|u_{0}|_{0}))t}
\]
uniformly in $N$. In particular,
\[
\|u(t)\|_{L^\infty}\leq |u_{0}|_{0}e^{-(1-d^{2}P_{1}(|u_{0}|_{0}))t}.
\]

\subsection{Uniqueness}
Due to the regularity of the solution constructed, the uniqueness follows from a standard contradiction argument. As this part is classical, we only sketch its proof. We consider two different solutions emanating from the same initial data, $u^{(1)}$ and $u^{(2)}$, and define
$$
v=u^{(1)}-u^{(2)}.
$$
Then
\[
\frac{1}{2}\frac{d}{dt}\|v\|_{L^2}^2+\|v_{,ij}\|_{L^2}^2=\int_{\mathbb{T}^d}\left(\frac{u^{(1)}_{,i} u^{(1)}_{,j}}{1+u^{(1)}}-\frac{u^{(2)}_{,i} u^{(2)}_{,j}}{1+u^{(2)}}\right)v_{,ij}dx.
\]
Integrating by parts and using H\"{o}lder inequality together with interpolation in Sobolev spaces we find the inequality
\[
\frac{d}{dt}\|v\|_{L^2}^2\leq c(\|u^{(1)}\|_{C^2}+\|u^{(2)}\|_{C^2})\|v\|_{L^2}^2.
\]
As the solution that we constructed is of class $L^2_TC^2$, we conclude the uniqueness of the solution.

\subsection{Spatial analyticity}
Let $\sigma\in (0,1)$ be a sufficiently small constant which will be fixed later. We define the function
\[
V^N(t,x)=e^{\sigma t\sqrt{-\Delta}}u^{N}(t,x) \quad \text{or equivalently} \quad \widehat{V^N}(t,k)=e^{\sigma t|k|}\widehat{u^N}(t,k).
\]
Let 
\eqn \label{Gu}
G(u)=\mathcal{P}_N\left(\left(u_{,i}u_{,j}\sum^{N}_{k=0}(-1)^{k}u^{k}\right)_{,ij}\right)
\een
Then, 
\eqn  \label{k Ud}
\partial_{t}\widehat{V^N}(t,k)&= -|k|^{4}\widehat{V^N}(t,k)-e^{\sigma t|k|}\widehat{G(u)}(t,k) +\sigma|k|\widehat{V^N}(t,k)
\een 
with $\widehat{V^N}(0,k)=\widehat{P_Nu_{0}}(k)$. We note that  
\eqn \label{product analyticity}
\left|e^{\sigma s\sqrt{-\Delta}}(fg)\right|_{0} \leq \sum_{k\in \mathbb{Z}^{d}} \sum_{l\in \mathbb{Z}^{d}}e^{\sigma s|k-l|} \left|\widehat{f}(k-l)\right|e^{\sigma s|l|}\left|\widehat{g}(l)\right| \leq \left|e^{\sigma s\sqrt{-\Delta}}f\right|_{0}\left|e^{\sigma s\sqrt{-\Delta}}g\right|_{0}.
\een 
In particular,
$$
\left|e^{\sigma s\sqrt{-\Delta}}(f^2)\right|_{0} \leq \sum_{k\in \mathbb{Z}^{d}} \sum_{l\in \mathbb{Z}^{d}}e^{\sigma s|k-l|} \left|\widehat{f}(k-l)\right|e^{\sigma s|l|}\left|\widehat{f}(l)\right| \leq \left|e^{\sigma s\sqrt{-\Delta}}f\right|_{0}^2.
$$
Furthermore,
\begin{align*}
\left|e^{\sigma s\sqrt{-\Delta}}(f^3)\right|_{0} &\leq \sum_{k\in \mathbb{Z}^{d}} \sum_{l\in \mathbb{Z}^{d}}e^{\sigma s|k-l|} \left|\widehat{f^{2}}(k-l)\right|e^{\sigma s|l|}\left|\widehat{f}(l)\right|\\
& \leq \left|e^{\sigma s\sqrt{-\Delta}}(f^2)\right|_{0}\left|e^{\sigma s\sqrt{-\Delta}}f\right|_{0}\\
& \leq \left|e^{\sigma s\sqrt{-\Delta}}f\right|_{0}^3.
\end{align*}
Similarly, 
\begin{align*}
\left|e^{\sigma s\sqrt{-\Delta}}(f^n)\right|_{0} \leq \left|e^{\sigma s\sqrt{-\Delta}}f\right|_{0}^n.
\end{align*}
We apply the previous estimates to the nonlinear terms in (\ref{k Ud}). Following the computations in Section \ref{apriori} used to obtain (\ref{DLSS final}), we derive
\begin{align*}
\left|e^{\sigma s\sqrt{-\Delta}}(u_{,i}u_{,j}u_{,ij})\right|_{0} \leq \left|e^{\sigma s\sqrt{-\Delta}}u\right|_{1}^2\left|e^{\sigma s\sqrt{-\Delta}}u\right|_{2}\leq\left|e^{\sigma s\sqrt{-\Delta}}u\right|_{0}^2\left|e^{\sigma s\sqrt{-\Delta}}u\right|_{4}.
\end{align*}
Performing the same estimates to the rest of nonlinear terms, we find the following inequality: 
\[
\frac{d}{dt}|V^N(t)|_{0}+|V^N(t)|_{4}\leq d^{2}P_{1}(|V^N(t)|_{0})|V^N(t)|_{4}+ \sigma |V^N(t)|_{1}.
\] 
Since $\langle V^N\rangle=\widehat{u^N}(0)=0$, we have $|V^N(t)|_{1}\leq |V^N(t)|_{4}$ and thus
\[
\frac{d}{dt}|V^N(t)|_{0}+(1-\sigma)|V^N(t)|_{4}\leq d^{2}P_{1}(|V^N(t)|_{0})|V^N(t)|_{4}.
\]
We now take $\sigma>0$ sufficiently small to satisfy $d^{2}P_{1}(|u_{0}|_{0})<1-\sigma$. Then, we obtain 
\eqn \label{analyticity final}
|V^N(t)|_{0}\leq |u_{0}|_{0} \quad \text{for all $t>0$.}
\een
Using \eqref{dualitywiener}, we conclude that 
$$
V^N\overset{\ast}{\rightharpoonup} V\in L^\infty_T\mathcal{A}_0.
$$
Now we observe that we can pass to the weak-$*$ limit in $N$ and conclude the desired bound for the function $V$. The convergence of $u^N$ implies that 
$$
V(t,x)=e^{\sigma t\sqrt{-\Delta}}u(t,x).
$$ 

%\begin{remark}
%In Section \ref{A: alternative analyticity}, we will provide another proof of analyticity of $u$ using the so-called Chemin-Lerner spaces. This, however, requires a slightly stronger condition on the size of $u_{0}$.
%\end{remark}

%%%%%%%%%%%%%%%%%%%%%%%%%%%%%%%%%%%%%%%%
\section{Proof of Theorem \ref{DLSS GW 2}} \label{sec 5}
%%%%%%%%%%%%%%%%%%%%%%%%%%%%%%%%%%%%%%%%
For the sake of notational simplicity, we write $\mu \Gamma=\gamma>0$ and $ \epsilon(2\Gamma^{2}-2\Gamma)=-\nu<0$ in (\ref{u ext}). We first rewrite (\ref{u ext}) as 
$$
u_{t}-\frac{\gamma}{3} u_{xxx}-\gamma ((u+1)(\log (u+1))_{xx})_{x}=-\nu((u+1)(\log (u+1))_{xx})_{xx}.
$$

%%%%%%%%%%%%%%%%%%%%%%
\subsection{Existence of solution}
%%%%%%%%%%%%%%%%%%%%%
Since $((u+1)(\log (u+1))_{xx})_{x}$ is a lower order term,  the first part of Theorem \ref{DLSS GW 2} is very similar to the one in Theorem \ref{DLSS GW 1} and so we only provide a priori estimates for the existence of solutions. Since  
\begin{equation*}
\begin{split}
& \left((1+u)(\log (u+1))_{xx}\right)_x=u_{xxx}-\frac{2u_{x}u_{xx}}{1+u}+\frac{u_{x}^3}{(1+u)^2},\\
& \left((1+u)(\log (u+1))_{xx}\right)_{xx}=u_{xxxx}-\frac{2(u_{x}u_{xxx}+u_{xx}^2)}{1+u}+\frac{2u_{x}^2u_{xx}+3u_{x}^2u_{xx}}{(1+u)^2}-\frac{2u_{x}^4}{(1+u)^3},
\end{split}
\end{equation*}
we take the Taylor expansion of the rational functions to derive  
\[
u_{t}-\frac{4\gamma}{3} u_{xxx}+\nu u_{xxxx}=\nu(\text{I}_1+\text{I}_2+\text{I}_3)+\gamma (\text{I}_4+\text{I}_5),
\]
with
\begin{align*}
\text{I}_1&=\frac{2u_{x}u_{xxx}+2u_{xx}^2}{1+u}=(2u_{x}u_{xxx}+2u_{xx}^2)\sum^{\infty}_{n=0}(-1)^{n}u^{n},\\
\text{I}_2&=-\frac{5u_{x}^2u_{xx}}{(1+u)^2}=5u_{x}^2u_{xx}\sum^{\infty}_{n=1}(-1)^{n}nu^{n-1},\\
\text{I}_3&=2\frac{u_{x}^4}{(1+u)^3}=2u_{x}^4\sum^{\infty}_{n=2}(-1)^{n}n(n-1)u^{n-2},\\
\text{I}_4&=-\frac{2u_{x}u_{xx}}{1+u}=-2u_{x}u_{xx}\sum^{\infty}_{n=0}(-1)^{n}u^{n},\\ \text{I}_5&=\frac{u_{x}^3}{(1+u)^2}=-u_{x}^3\sum^{\infty}_{n=1}(-1)^{n}nu^{n-1}.
\end{align*}
Since 
\eqn \label{vanishing 3}
\frac{Re(\overline{\widehat{u}(t,k)}\widehat{u_{xxx}}(t,k))}{|\widehat{u}(t,k)|}=0, 
\een
we have  
\[
\frac{d}{dt}|u(t)|_{0}+\nu |u(t)|_{4}\leq \nu \sum^{3}_{j=1}|\text{I}_j|_{0}+\gamma \sum^{5}_{j=4}|\text{I}_{j}|_{0}. 
\]
We recall $P_{1}$ and $P_{2}$ in (\ref{P1}) and (\ref{P2}):
\[
P_{1}(z)=\frac{4z}{1- z}+\frac{5z^2}{(1- z)^{2}}+\frac{4z^{3}}{(1- z)^{3}}, \quad P_2(z)=\frac{z}{1- z}+\frac{z^2}{(1- z)^2}.
\]
By applying (\ref{interpolation}) and \eqref{geometric series}, we obtain 
\begin{equation*}
\begin{split}
& \sum^{3}_{i=1}\left|\text{I}_i\right|_{0} \leq \left(\frac{4|u(t)|_{0}}{1- |u|_{0}}+\frac{5|u(t)|_{0}^2}{(1- |u|_{0})^2}+\frac{4|u(t)|_{0}^3}{(1- |u|_{0})^3}\right)|u|_{4}=P_1(|u|_{0})|u|_{4},\\
& \sum^{5}_{i=4}\left|\text{I}_i\right|_{0}\leq \left(\frac{|u(t)|_{0}}{1- |u|_{0}}+\frac{|u(t)|_{0}^2}{(1- |u|_{0})^2}\right)|u|_{3}\leq P_2(|u|_{0})|u|_{4}
\end{split}
\end{equation*}
which imply that
\eqn \label{eODE}
\frac{d}{dt}|u(t)|_{0}+\nu |u(t)|_{4}\leq \nu P_1(|u|_{0})|u|_{4}+ \gamma P_2(|u|_{0})|u|_{4}. 
\een
Suppose $u_{0}$ satisfies 
\[
|u_{0}|_{0}<1 \quad \text{and} \quad \gamma P_{2}(|u_{0}|_{0})+\nu P_{1}(|u_{0}|_{0})<\nu.
\]
Then, we have 
\[
|u(t)|_{0} + \left(\nu-\gamma P_{2}(|u_{0}|_{0})-\nu P_{1}(|u_{0}|_{0})\right) \int^{t}_{0}|u(s)|_{4}ds \leq |u_{0}|_{0} \quad \text{for all $t\ge 0$.}
\]

%%%%%%%%%%%%%%%%%%%%%%%%%
\subsection{Decay to equilibrium}
%%%%%%%%%%%%%%%%%%%%%%%%%
As before, since $\langle u(t) \rangle=0$, we have $|u(t)|_{0}\leq |u(t)|_{4}$. By (\ref{eODE}),  
\[
\|u(t)\|_{L^\infty}\leq |u_{0}|_{0}e^{-(\nu-\gamma P_{2}(|u_{0}|_{0})-\nu P_{1}(|u_{0}|_{0}))t}.
\]

\subsection{Uniqueness}
The uniqueness follows similarly as before.

%%%%%%%%%%%%%%%%%%%%%%%%%%%%
\subsection{Spatial analyticity}
%%%%%%%%%%%%%%%%%%%%%%%%%%%%
Let $\sigma\in (0,1)$ be a small parameter that will be fixed later. As the approximation process mimics the method of constructing approximate solutions in Theorem \ref{DLSS GW 1}, we skip the approximation process. Thus, we continue directly to obtain appropriate estimates for 
$$
V(t,x)=e^{\sigma t\sqrt{-\Delta}}u(t,x).
$$ 
We observe that $V$ satisfies 
\[ 
\partial_{t}\widehat{V}(t,k) = -|k|^{4}\widehat{V}(t,k) -\frac{4\gamma}{3}\widehat{V_{xxx}}(t,k)+\nu e^{\sigma t|k|}(\text{I}_1+\text{I}_2+\text{I}_3)+\gamma e^{\sigma t|k|}(\text{I}_4+\text{I}_5)+\sigma|k|\widehat{V}(t,k)
\] 
with $\widehat{V}(0,k)=\widehat{u_{0}}(k)$. Using (\ref{vanishing 3}) applied to $V$, $|V(t)|_{1}\leq |V(t)|_{4}$, and (\ref{product analyticity}), we obtain 
\[
\frac{d}{dt}|V(t)|_{0}+\nu |V(t)|_{4}\leq \left(\gamma P_{3}(|V(t)|_{0})+\nu P_{4}(|V(t)|_{0})\right)|V(t)|_{4} +\sigma |V(t)|_{4}.
\]
Suppose $u_{0}$ is such that 
$$
\nu-\sigma-\gamma P_{2}(|u_{0}|_{0})-\nu P_{1}(|u_{0}|_{0})>0.
$$
Then, we have 
\eqn \label{gDLSS uniform bound}
|V(t)|_{0} \leq |u_{0}|_{0} \quad \text{for all $t>0$.}
\een
This completes the proof of Theorem \ref{DLSS GW 2}.

%%%%%%%%%%%%%%%%%%%
\section*{Acknowledgments}
%%%%%%%%%%%%%%%%%%%
H. Bae was supported by NRF-2018R1D1A1B07049015.   Part of this research was conducted while R. Granero-Belinch\'on was visiting UNIST funded by NRF-2018R1D1A1B07049015.

\bibliographystyle{plain}
%\bibliography{references}

\end{document}